\documentclass[11pt]{amsart}
\usepackage{amssymb}
\usepackage{amsfonts}
\usepackage{amscd}
\usepackage{amsmath}
\usepackage{mathrsfs}
\usepackage{curves}
\usepackage{epsfig}
\usepackage[pass]{geometry}
\usepackage{graphicx}
\usepackage{verbatim}
\usepackage{latexsym}
\usepackage{eucal}
\usepackage{bbm}
\usepackage{tikz}
\usetikzlibrary{matrix}

\newtheorem{theorem}{Theorem}[section]
\newtheorem{lemma}[theorem]{Lemma}

\newtheorem{cor}[theorem]{Corollary}

\newtheorem{definition}[theorem]{Definition}

\newtheorem{remark}[theorem]{Remark}

\def \mcd {{\mathscr D}}
\def \mce {{\mathscr E}}

\def \mcx {{\mathscr X}}
\def \mcy {{\mathscr Y}}
\def \mcz {{\mathscr Z}}

\def \mbc {{\mathbb C}}

\def \mbn {{\mathbb N}}
\def \mbr {{\mathbb R}}

\def \re {\operatorname{Re}}

\def \defeq {\stackrel{\operatorname{def}}{=}}
\def \beqq {\begin{equation}}
\def \eeqq {\end{equation}}
\def \WF {\text{WF}}

\def \bpf {\begin{proof}}
\def \epf {\end{proof}}
\def \beq {\begin{equation*}}
\def \eeq {\end{equation*}}

\def \eps {\epsilon}   
\def \la {\lambda}   
\def \La {\Lambda}    
   
\def \lap {\Delta}
\def \p {\partial}
\def \ha {\frac{1}{2}}

\def\mr {\mbr}

\def \CI {{C^\infty}}

\def \ido {{\stackrel{o}{\operatorname{I}}}}

\begin{document}
\title[]{Identification of nonlinear  beam-hardening effects in X-ray tomography}
\author[]{Yiran Wang}
\address{Yiran Wang
\newline
\indent Department of Mathematics, Emory University
\newline
\indent 400 Dowman Drive, Atlanta, Georgia 30322, USA}
\email{yiran.wang@emory.edu}
\begin{abstract}
We study streaking artifacts caused by beam-hardening effects in X-ray computed tomography (CT). The effect is known to be nonlinear. We show that the nonlinearity can be recovered from the observed artifacts for strictly convex bodies. The result provides a theoretical support for removal of the artifacts. 
\end{abstract}
\date{\today}
 
\maketitle

\section{Introduction}
In X-ray computed tomography (CT), artifacts due to beam-hardening effects is common for patients with medical implants. They are notorious for  causing degradation of CT images and difficulties for diagnosis. The reduction or removal of such artifacts has drawn numerous research efforts, but the problem still remains  one of the major challenges in X-ray CT. 

\begin{figure}[t]
\centering 
\includegraphics[scale=.52]{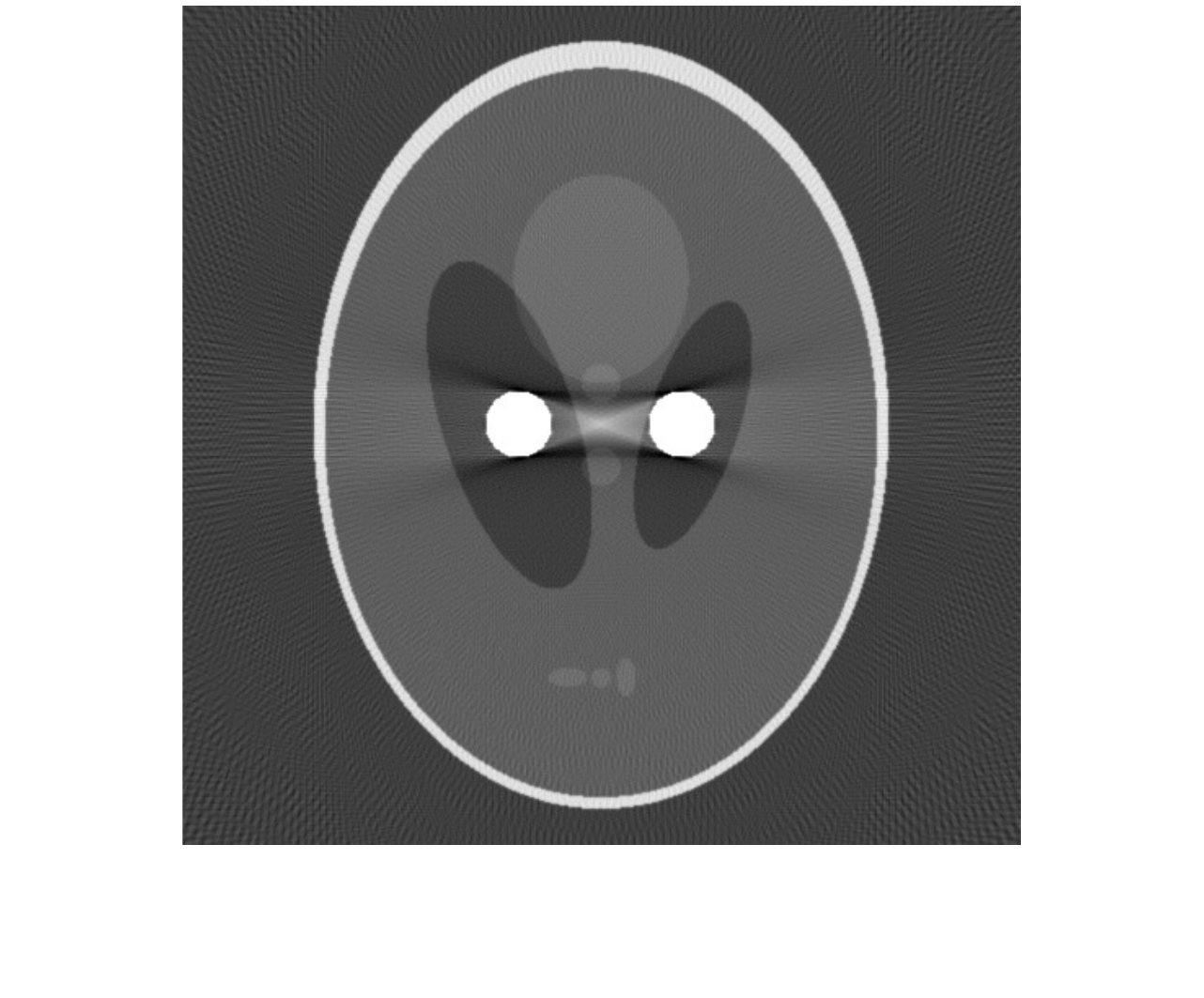}   \includegraphics[scale=0.52]{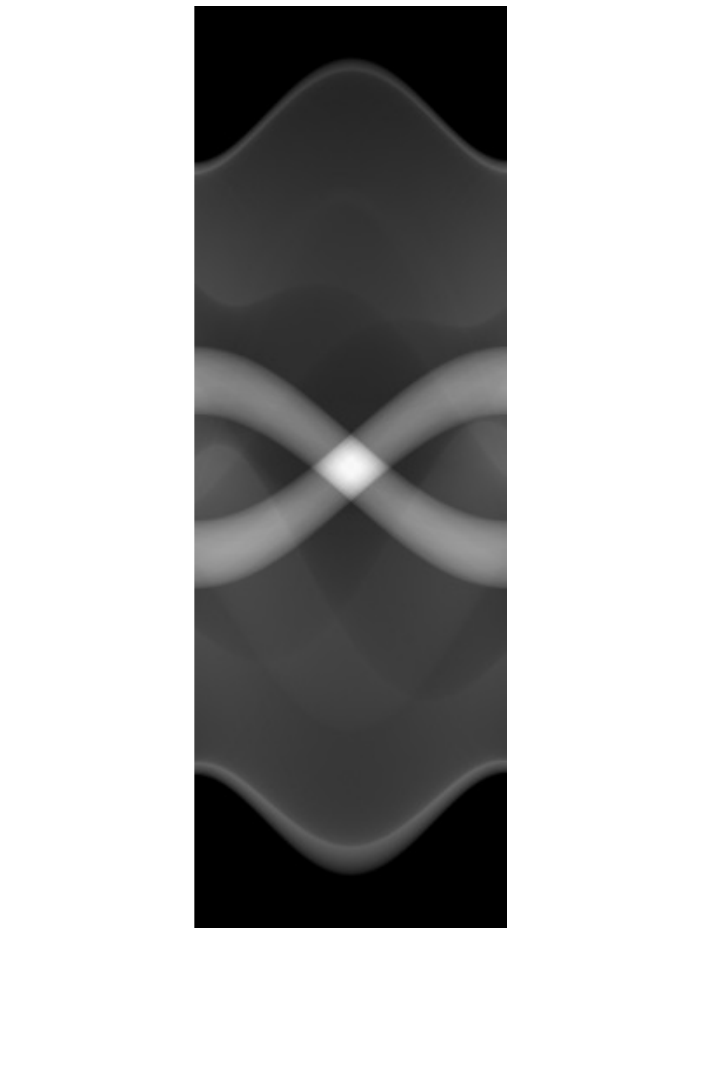}
\vspace{-0.7cm}
\caption{Illustration of the beam-hardening artifacts produced by a quadratic nonlinearity. The left figure shows the reconstruction from FBP. The two bright disks represents the metal objects. The artifacts appear to be straight lines tangent to both disks. The right figure shows the sinogram (color map of $P$ in \eqref{eq-ctp}).  The two bright strips correspond to the support of $R\chi_D.$ }
\label{fig-ex}
\end{figure}

In  \cite{Seo}, the authors demonstrated the nonlinear nature of the beam-hardening effects. Let $f$ be the attenuation coefficient of the object being imaged. Because of the polychromatic nature of X-ray beams, we consider the dependency of $f$ on the energy level $E$. This is particularly significant for metal objects. Assume that   $E \in [E_0 -\eps, E_0 + \eps]$, $\eps > 0$, and write $f$ as $f_E$. We write  $f_E(x) = f_{E_0}(x) +  \alpha \chi_D$, where $\chi_D$ is the characteristic function for a metal object $D\subset \mbr^2$, and $\alpha > 0$ a constant which can be thought of as the approximation of the derivative of $f_E$ in $E$. The X-ray data can be derived from the Beer-Lambert law which gives  
\beqq\label{eq-ctp}
P = R f_{E} + P_{MA}. 
\eeqq 
Here,  $Rf_{E}$ denote the Radon transform of $f_E$, and $P_{MA}$ denotes a mismatch term. Under further assumptions, it is derived in \cite{Seo} that 
\beqq\label{eq-ma}
P_{MA} = - \ln\left(\frac{\sinh(\alpha \eps R \chi_D)}{\alpha \eps R \chi_D}\right)
\eeqq
is a nonlinear function of $R\chi_D$. If one applies the filtered back-projection (FBP) to reconstruct $f_E$,  the mismatch term $P_{MA}$ leads to the artifact, see Figure \ref{fig-ex} for a demonstration.  By using the notion of wave front set in microlocal analysis, the authors of \cite{Seo} gave a mathematical characterization of the artifacts. For strictly convex objects, the artifacts appear to be straight lines tangent to at least two boundary points of the metal objects, see Figure \ref{fig-ex}. For more complicated situations, the artifacts and their relation to the geometry of metal regions are further studied in \cite{PUW, WaZo}. Finally, we mention that artifacts due to similar  mechanism are also known in the attenuated X-ray tomography, see \cite{Kat}.

The identification of the nonlinear effect in $P_{MA}$ is the key for removal of  the artifacts. In practice, the shape of the metal objects can usually be acquired so it is reasonable to assume that $\chi_D$ is known and think of $P_{MA}$ as a nonlinear function $F(R{\chi_D})$. Then one can remove $P_{MA}$ if $F$ is known. For example, the numerical scheme  developed in \cite{Seo1} consists of two steps: first, one recovers $\chi_D$ from  the reconstruction $f_{CT}$ using image segmentation techniques; second, one can use the model \eqref{eq-ma} and find the optimal $\alpha$ which reduces the artifact. In general, not much can be said about the nonlinear function as it depends on many factors such as the energy distribution of X-ray beams, geometry and physical properties of metal objects and even $f_{E_0}$. Lately, there have been increasing efforts to find the nonlinear effect by using deep learning techniques, see for example \cite{Seo2, ZhYu}.  

The purpose of this note is to show that the nonlinearity can be identified from the observed artifacts. Also, we provide a constructive proof which is  potentially useful in practice.

\section{The main result}\label{sec-main}
Because of the local nature of the problem as we explain later, it suffices to state our main result in a relatively simple setting. We assume that the metal region $D = D_1\cup D_2$,  where $D_j, j=1, 2$ are simply connected disjoint bounded domains in $\mbr^2$ with smooth boundary $\p D_j$. 
Let $\chi_D =  \chi_{D_1} + \chi_{D_2}$ be the characteristic function of $D$. We assume that the attenuation coefficient is of the form  
\beqq\label{eq-f}
f(x) = h(x) + \chi_D(x), \quad x\in \mbr^2
\eeqq
where $h$ represents the attenuation coefficient of the tissue. 
For the Radon transform on $\mbr^2$, we use the following parametrization 
\[
R f(s, \phi) = \int_{x_1\cos\phi + x_2\sin\phi = s} f(x) dx 
\]
Here, $x\in \mbr$ and $(s, \phi)\in M = \mbr \times (-\pi, \pi)$. Note that  because we are dealing with a local problem for $Rf$, it suffices to use local coordinates in the region of interest. We model the beam-hardening effects by a polynomial function $F: \mbc\rightarrow \mbc$ of the form 
\beqq\label{eq-F}
F(t) = \sum_{j = 2}^J a_j t^j
\eeqq
where $a_j$ are constant and $J\in \mbn.$ 
Then the X-ray CT data is modeled by 
 \beqq\label{eq-pct}
P = R f + P_{MA}, \quad P_{MA} =  F(R \chi_D)
\eeqq
For the reconstruction, we apply the FBP to get 
\beqq\label{eq-fct}
f_{CT} =   f  +  f_{MA}, \ \ f_{MA} = R^*I^{-1} P_{MA}
 \eeqq
where  $I^{-1}$ is the Riesz potential and $R^*$ denotes the adjoint of $R$, see for example \cite{Qui}. 
Our main result is 
\begin{theorem}\label{thm-main}
Suppose $D_j, j = 1, 2$ are strictly convex and $h\in C_0^\infty(\mbr^2)$. Then $F$ in \eqref{eq-F} is uniquely determined by $f_{CT}$ in \eqref{eq-fct}. 
\end{theorem}

In fact, as we will see in the proof of Theorem \ref{thm-main}, it suffices to use $f_{CT}$ away from $\p D$ to determine $F$. So the theorem really says that $F$ can be determined from the artifacts. This is important because $f$ in \eqref{eq-f} is singular at $\p D$. The singularities of $f_{MA}$ at $\p D$ are weaker so it would be more difficult to recover $F$ from those singularities. Among other things, one useful consequence is the following
\begin{cor}\label{cor}
Under the assumption of Theorem \ref{thm-main}, $f_{CT}\in C^\infty(\mbr^2\backslash \p D)$ if and only if $F = 0.$
\end{cor}
The result implies that the artifact is always visible unless there is no beam-hardening effect. This seems to be the first existence result for the streaking artifacts. 

Next, we make a few remarks regarding the assumptions of the results. 
\begin{remark}[The nonlinear function $F$]
As studied in \cite{Seo, PUW, WaZo}, the artifacts are generated from $P_{MA}$ in \eqref{eq-pct} at certain discrete points. We will see in Section \ref{sec-metal} that near those points, it suffices to consider $t$ small in \eqref{eq-F}. Thus assuming $F$ of the form \eqref{eq-F} is not  restrictive. In fact, Theorem \ref{thm-main} can be applied to any of those points to reconstruct the coefficients of the Taylor series expansion of  a general nonlinear function.   Note that any linear term in $F$ can be thought of as part of $f$ in \eqref{eq-f}, that's why they are not included in \eqref{eq-F}. 
\end{remark}

\begin{remark}[The number of objects]
The streaking artifacts are  associated with lines tangent to both $D_1$ and $D_2$. If the metal regions consists of more than two disjoint simply connected regions, our results apply as long as there is no line tangent to more than two of the regions. We refer to \cite{PUW} for the detailed treatment of that case. 
\end{remark}

\begin{remark}[The convexity assumption]
The assumption that  $D_j, j = 1, 2$ are strictly convex is essential. In Section \ref{sec-radon}, we will show that $R\chi_{D_j}, j = 1, 2$ possess better regularity properties under the strict convexity assumption which is the key to determine $F$. In fact, if $D_j$ are not strictly convex, it seems impossible to recover $F$ from the artifacts. 
\end{remark}

\begin{remark}[The regularity of $h$] 
We assumed $h\in C_0^\infty$ for simplicity. Actually, the results hold as long as the singular support of $h$ does not meet the streaking artifacts. 
\end{remark} 
 
Finally, we briefly discuss the ideas of the proof. We already mentioned that  streaking artifacts are associated with singularities, more precisely wave front sets of  $F(R\chi_D)$  due to  the nonlinear interactions of the singularities in $R\chi_D.$ With more precise notions of  Lagrangian distributions, quantitative results on the strength of the artifacts are obtained in \cite{PUW}, see also \cite{WaZo} for non-convex objects. Essentially, these results provided the upper bound of the wave front set  which tells where the artifact {\em could} appear.  
Our  idea is that when the artifact {\em actually} happens, it carries information of the nonlinear function $F$. So we can use the artifact to reconstruct $F$. The philosophy that nonlinearity can help solving inverse problems was perhaps first demonstrated in \cite{KLU} for nonlinear wave equations. In the last decade, the method has undergone rapid developments and the majority of the work relies on the idea of higher order linearization. Unfortunately, this method is not applicable to our problem because it requires the X-ray data for a family of metal objects but we only have one. 

The way we overcome the  difficulty is to study the fine structure of  singularities in $F(R\chi_D)$. A key observation, discussed in Section \ref{sec-radon} is that when $D$ is strictly convex, $R\chi_D$ can be locally written as an asymptotic summation of conormal distributions with increasing regularities. This allows us to obtain expansion of $F(R\chi_D)$ near the interaction point in terms of the strength of singularities instead of the magnitudes in the higher order linearization expansion. To recover $F$, another key component is to show that all terms in the expansion are non-trivial. This   is done in Section \ref{sec-non}, where we finish the proof of Theorem \ref{thm-main} and Corollary \ref{cor}.

\section{Microlocal analysis of the artifacts}\label{sec-metal}
In this section, we consider the artifact generation for a quadratic nonlinearity. Some of the analysis appeared in \cite{PUW} but we need to improve them so they can be used for more general nonlinearities. 
 
\subsection{Regularity of $R\chi_D$.}\label{sec-reg} We start with the notion of conormal distributions, see Section 18.2 of \cite{Ho3} for details. Let  $\Omega \subset \mr^{n}, n\in \mbn$  denote an open and relatively compact subset. Let  $\Sigma \subset \Omega$ be a submanifold of co-dimension $k$. The set of conormal distributions of order $m$ is denoted by $I^{m}(\Omega, \Sigma)$. Such distributions can be defined in a coordinate invariant way but we only need the local representations. 
%
 According to Theorem 18.2.8 of \cite{Ho3}, $u\in I^{m}(\Omega, \Sigma)$ if and only if  $u \in C^\infty(\mr^{n}\setminus \Sigma)$ and near any point $p\in \Sigma$ and in local coordinates where 
  $\Sigma=\{y_1=y_2=\ldots= y_k=0\},$  $y=(y', y''),$  $y'=(y_1, y_2, \ldots, y_k),$ $y''\in \mr^{n-k},$ we have
\begin{gather}
u(y)= \int_{\mr^k} e^{i  y'\cdot \eta' } a(\eta',y'') \; d\eta', \;\ a\in S^{m+\frac{n-2k}{4}}( \mr^{k}  \times \mr^{n-k} ), \label{con1}
\end{gather}
where for $r\in \mr,$ $S^{r}(\mr^k \times \mr^{n-k})$ is the class of symbols satisfying
\begin{gather*}
|\p_{y''}^\alpha \p_{\eta'}^\beta a(\eta',y'')| \leq C_{\alpha, \beta}(1+|\eta'|)^{r-|\beta|}.
\end{gather*}
It is known that $\WF(u)\subset N^*\Sigma$ and the space of distributions satisfy
\begin{gather*}
I^{m}(\Omega, \Sigma) \subset I^{m'}(\Omega, \Sigma), \;\ m<m'.
\end{gather*}
The principal symbol of $u$ is defined to be  the equivalence class of $a(\eta', y'')$ in the quotient $S^{m+\frac{n-2k}{4}}(\mr^k\times \mr^{n-k})/ S^{m+\frac{n-2k}{4}-1}(\mr^k \times \mr^{n-k})$ and the map
\begin{gather*}
 I^{m}(\Omega, \Sigma)/ I^{m-1}(\Omega, \Sigma) \longrightarrow S^{m+\frac{n-2k}{4}}(\mr^k\times \mr^{n-k})/ S^{m+\frac{n-2k}{4}-1}(\mr^k \times \mr^{n-k})\\
 [u] \longmapsto [a],
 \end{gather*}
is an isomorphism. The symbol map can be invariantly defined as in \cite{Ho3}, but since our analysis is completely local, we do not need to discuss that. Below, we mostly  consider $n = 2$ and $k=1.$ 

 Let $D_j, j = 1, 2$ be a simply connected bounded domain in $\mbr^2$ with smooth strictly convex boundary $\p D_j$. The characteristic function $\chi_{D_j} \in I^{-1}(\mbr^2;  \p D_j)$. 
It is known that $R \chi_{D_j}$ is a conormal distribution, see for example \cite{PUW}. We will see a direct calculation in Section 3. To describe the conormal distribution, we start with 
 $R: \mce'(\mbr^2)\rightarrow \mcd'(M)$ as an elliptic Fourier integral operator. Using local coordinates $(s, \phi)$ for $M$ and $x = (x_1, x_2)$ for $\mbr^2$, we write the Schwartz kernel of $R$, denoted by $K_R$, as an oscillatory integral
\beq
K_R(s, \phi, x) = \frac{1}{(2\pi)^\ha} \int_\mbr  e^{i(  x_1 \cos \phi + x_2 \sin \phi -s)\la} d\la.
\eeq
The phase function is $\phi(s, \theta, x; \la) = (x_1\cos\phi + x_2\sin \phi -s)\la$ so the associated Lagrangian submanifold of $T^*M\times T^*\mbr^2$ is 
\beq
\begin{gathered}
\La  = \{(x_1 \cos \phi + x_2 \sin \phi, \phi, -\la,  \la(-x_1\sin \phi + x_2 \cos \phi); x_1, x_2,  \la \cos \phi, \la \sin \phi):\\
 \la \in \mbr\backslash 0,  \phi \in (-\pi, \pi), x_1, x_2\in \mbr\}.
 \end{gathered}
\eeq
In particular, $K_R \in I^{-\ha}(M\times\mbr^2; \La)$. We denote the homogeneous canonical relation by 
\beqq\label{eqcaC}
\begin{gathered}
C  
=  \{(x_1 \cos \phi + x_2 \sin \phi, \phi, -\la,  \la(-x_1\sin \phi + x_2 \cos \phi); \\
x_1, x_2,  -\la \cos \phi, - \la \sin \phi):\\
 \la \in \mbr\backslash 0,  \phi \in (-\pi, \pi), x_1, x_2\in \mbr\} \subset T^*M \backslash 0 \times T^*\mbr^2\backslash 0
\end{gathered}
\eeqq 
Let $C_j \defeq  N^*\p D_j\backslash 0, j = 1, 2$ and we think of them as canonical relations of $\chi_{D_j}$, see Appendix \ref{sec-app}. The composition of the two homogeneous canonical relations $C, C_j$ is transversal (see Appendix \ref{sec-app}) so the composition $C\circ C_j$ is a homogeneous canonical relation which is   a Lagrangian submanifold of $T^*M$.  Under the strict convexity assumption, the Lagrangian becomes a conormal bundle. In fact, the projection of $C\circ C_j$ to $M$ is injective and the projection is 
\beqq\label{eq-sj}
S_j  \defeq \{(s, \phi) \in M : s = x_1\cos \phi + x_2\sin\phi, \ \  (x_1, x_2) \in \p D_j \}
\eeqq
which is a co-dimension one submanifolds of $M$. We have $C\circ C_j =  N^*S_j\backslash 0, j = 1, 2$.  One can apply the FIO composition theorem \cite[Theorem 25.2.3]{Ho4} to conclude that $R \chi_{D_j} \in I^{-\frac 32}(M;  S_j), j = 1, 2.$

\subsection{The nonlinear analysis.} We consider a quadratic nonlinearity $F(t) = t^2$ in \eqref{eq-F} and let 
\beqq\label{eq-rd}
\tilde P_{MA}  \defeq   (R(\chi_D))^2 =  (R(\chi_{D_1}))^2 + (R(\chi_{D_2}))^2 + 2 R(\chi_{D_1})R(\chi_{D_2}).
\eeqq
We analyze the singularity  in each terms. For $(R \chi_{D_j})^2, j = 1, 2$, we recall the following multiplicative property of conormal distributions, see   \cite{Pir}: 
\begin{lemma}\label{lm-multi} 
If $\Sigma\subset \Omega$ is a $C^\infty$ hypersurface, if $u, v\in I^{m-\frac{n}{4}+\ha}(\Omega,\Sigma)$ and $m<-1,$ then $uv\in I^{m-\frac{n}{4}+\ha}(\Omega,\Sigma).$
\end{lemma}
We conclude that $(R \chi_{D_j})^2 \in I^{-\frac 32}(M;  S_j), j = 1, 2$. So $\WF((R \chi_{D_j})^2) \subset N^*S_j$ does not produce new singularities.

Next consider the product $R\chi_{D_1}R\chi_{D_2}$ in \eqref{eq-rd}. This term will produce new singularities and the result can be described by using the notion of paired Lagrangian distributions, see \cite{GU} for details. Let $X$ be an $n$-dimensional manifold, $\La_0, \La_1\subset T^*X\backslash 0$ be two cleanly intersecting Lagrangians in the sense that $\Sigma = \La_0\cap \La_1$ is smooth and $T_q\Sigma = T_q\La_0 \cap T_q\La_1, q\in \Sigma.$ The paired Lagrangian distributions associated with the pair $(\La_0, \La_1)$ with order $p, l\in \mbr$ is denoted by $I^{p, l}(\La_0, \La_1)$. We only need the case when $\La_0, \La_1$ are conormal bundles. Locally, such distributions can be described as follows, see \cite{GrU93}. Let $x = (x_1, \cdots, x_n)$ be  coordinates of $\mbr^n.$ Let $k_1, k_2 \in \mbn$ and $k_1 + k_2 \leq n.$ Consider 
\beq
\begin{gathered}
Y_1 = \{x_1= x_2 = \cdots x_{k_1} = 0\} = \{x' = 0\},\\
Y_2 = \{x_1 = x_2 = \cdots x_{k_1 + k_2 = 0}\} = \{x' = 0, x'' = 0\}.
\end{gathered}
\eeq
Let $\La_0 = N^*Y_1\backslash 0, \La_1 = N^*Y_2\backslash 0$ and $u\in I^{p, l}(\La_0, \La_1)$ can be written as 
\beqq\label{eq-pair}
u(x) = \int_{\mbr^{k_1 + k_2}} e^{i(x'\cdot \xi' + x'' \cdot \xi'')} a(x, \xi', \xi'') d\xi' d\xi''
\eeqq
with $a(x; \xi', \xi'')$ belonging to the product type symbols
\beq
\begin{gathered}
S^{\mu, \mu'}(\mbr^n\times (\mbr^{k_1} \backslash 0) \times \mbr^{k_2}) = \{ a\in C^\infty: |\p_x^\gamma \p_{\xi''}^\beta\p_{\xi'}^\alpha a(x, \xi)|\\
\leq C_{\alpha\beta\gamma K} \langle \xi', \xi''\rangle^{\mu - |\alpha|} \langle \xi''\rangle^{\mu' - |\beta|} \}
\end{gathered}
\eeq
where $\mu  = p-k_1/2 + n/4, \mu' = l - k_2/2$. 
We recall the fact that  if $u\in I^{p, l}(\La_0, \La_1)$, then $\WF(u)\subset \La_0\cup \La_1$. Also, $u\in I^{p+l}(\La_0\backslash \La_1)$ and $u\in I^{p}(\La_1\backslash \La_0)$ as Lagrangian distributions. Thus the principal symbols can be defined invariantly for each piece. In fact, one can also define the notion of principal symbols for $u$ invariantly, see \cite{GU}. However, we only need the behavior of these distributions locally. So it suffices to work with the expression \eqref{eq-pair}

To analyze the singularities in $R\chi_{D_1}R\chi_{D_2}$, we need to know how the Lagrangians intersects. As shown in \cite{PUW}, if $D_1, D_2$ are strictly convex, then  $S_1$ intersect $S_2$ transversally at a finite point set $S_\diamond = S_1\cap S_2$, see Figure \ref{fig-ex}.  
\begin{lemma}
For each $q\in S_\diamond$, there exists neighborhood $O$ of $q$ such that  in $O$, 
\beqq\label{lm-dij}
R(\chi_{D_1})R(\chi_{D_2}) \in  I^{-\frac 32, -1}(T_q^*M, N^*S_1) + I^{-\frac 32, -1}(T_q^*M, N^*S_2)
\eeqq
and the principle symbol on $T_q^*M \backslash (N^*S_1\cup N^*S_2)$  is non-vanishing.
\end{lemma}
\bpf
We repeat the proof of Lemma 1.1 of \cite{GrU93}. Consider the intersection of $S_1, S_2$ at $q.$ We choose local coordinates $(x_1, x_2)$ for $\mbr^2$ such that $q = (0, 0)$, $S_1 = \{x_1 = 0\}$ and $S_2 = \{x_2 =  0\}$. Then we can write 
\beq
\begin{gathered}
R\chi_{D_1}(x) = \int_{\mbr} e^{ix_1\xi_1} a(x, \xi_1) d\xi_1, \quad a\in S^{-3/2}(\mbr^2\times (\mbr\backslash0))\\
R\chi_{D_2}(x) = \int_{\mbr} e^{ix_2\xi_2} b(x, \xi_2) d\xi_2, \quad b\in S^{-3/2}(\mbr^2\times (\mbr\backslash0))
\end{gathered}
\eeq 
and the principal symbols of $a, b$ are non-vanishing. Then we get 
\beqq\label{eq-p1}
R\chi_{D_1}(x)R\chi_{D_2}(x) =  \int_{\mbr} e^{ix_1\xi_1 + i x_2\xi_2} a(x, \xi_1)b(x, \xi_2) d\xi_1d\xi_2
\eeqq
Introduce a cutoff function $\chi(t)\in C_0^\infty(\mbr)$, $\chi = 1$ for $|t|\leq 1/2$ and $\chi = 0$ for $|t| \geq 1$. Then we have 
\beq
\begin{gathered}
\chi(\langle \xi_2\rangle/\langle \xi_1\rangle) a(x, \xi_1)b(x, \xi_2) \in S^{-3/2, -3/2} (\mbr^2\times (\mbr\backslash 0)\times \mbr)\\
(1 -\chi)(\langle \xi_2\rangle/\langle \xi_1\rangle) a(x, \xi_1)b(x, \xi_2) \in S^{-3/2, -3/2} (\mbr^2\times (\mbr\backslash 0)\times \mbr)
\end{gathered}
\eeq
So the product is a sum of two Lagrangian distributions with orders $p = -3/2 + 1/2 - 1/2 = -3/2$ and $l = -3/2 + 1/2 = -1.$

To find the principal symbol on $T_q^*M \backslash (N^*S_1\cup N^*S_2)$, we consider \eqref{eq-p1} for $C_1\langle \xi_2\rangle 
\leq \langle \xi_1\rangle \leq C_2 \langle \xi_2\rangle$ for some positive constants $C_1, C_2$. Then the symbol $ab\in S^{-3}(\mbr^2\times \mbr^2)$ and the principal symbol is given by the product of principal symbols of $a$ and $b$.  
\epf

In conclusion, we proved that $\WF(\tilde P_{MA})\subset (\cup_{q\in S_\diamond} T_q^*M) \cup N^*S_1\cup N^*S_2$. Moreover, $(\cup_{q\in S_\diamond}T_q^*M) \backslash N^*(S_1\cup S_2) \subset \WF(\tilde P_{MA})$.

\subsection{Description of the artifact}\label{subsec-art}
Consider $\tilde f_{MA}  \defeq    R^* I^{-1} (R \chi_D)^{2 }.$ 
We show that this term contributes to the streaking artifacts. We define  $L_\diamond = \{L: \text{$L$ is a straight line in $\mbr^2$ tangent to $D_1$ and $D_2$}\}.$ 
 
\begin{lemma}\label{lm-fbp}
Fix any $L\in L_\diamond$, consider $\tilde f_{MA}$ near $L$ and away from $\p D_1\cup  \p D_2 \cup L_\diamond \backslash \{L\}$. Then we have  $\tilde f_{MA} \in I^{-2}(N^*L)$  and the principal symbol is non-vanishing. 
\end{lemma}
\bpf
Essentially, $I^{-1}$ is an pseudo-differential operator of order $1$, see e.g. \cite{WaZo}. Also, we know that $R^*$ is an elliptic FIO of order $-\ha$. Let $C^*$ be the canonical relation of $R^*.$  Then we check that $C^*\circ N^*S_{j}\backslash 0 = C^*\circ C\circ N^*\p D_j\backslash 0=  N^*\p D_j\backslash 0, j= 1, 2.$ It follows from the wave front analysis that 
\beq
\WF(R^*I^{-1}((R\chi_{D_j}))^2) \subset N^*\p D_j, \quad j = 1, 2.
\eeq

Next, consider singularities in $R^*I^{-1}(R\chi_{D_1}  R\chi_{D_2})$. First of all, for $q\in S_\diamond$ and  let $\chi_q$ be a smooth cut-off function supported near $q$, we have 
\beq
I^{-1}\big( \chi_q R(\chi_{D_1})R(\chi_{D_2}) \big)  \in I^{-\ha, -1}(T_q^*M, N^*S_1) + I^{-\ha, -1}(T_q^*M, N^*S_2) 
\eeq
and the principal symbol at $T_q^*M\backslash (N^*S_1\cup N^*S_2)$ is non-vanishing. Here, we used Proposition 4.1 of  \cite{GU}. For the application of $R^*$, we can still use Proposition 4.1 of \cite{GU}. The transversality of the compositions  $C^*\circ (T_q^*M\backslash 0)$ and $C^*\circ (N^*S_j\backslash 0), j = 1, 2$ are verified in Appendix \ref{sec-app}. In particular, 
$C^*\circ (T_q^*M\backslash 0)= N^*L\backslash 0$, 
where
\beq  
 L = \{x\in \mbr^2: x_1\cos\phi + x_2\sin \phi = s, q = (s, \phi)\}\in L_\diamond. 
\eeq  
So we get 
\beq
R^*\circ I^{-1}\big(\chi_q R(\chi_{D_1})R(\chi_{D_2})\big)\in I^{-1, -1}(N^*L, N^*\p D_1) + I^{-1, -1}(N^*L, N^*\p D_2).
\eeq
The principal symbol at $N^*L\backslash (N^*\p D_1\cup N^*D_2)$ is the product of the principal symbols of $R^*I^{-1}$ and $\chi_q R(\chi_{D_1})R(\chi_{D_2})$ at $N^*L$ so it is non-vanishing. The analysis can be repeated for each $q \in S_\diamond$ and the proof is completed. 
\epf

According to \cite[Definition 3.2]{Seo}, the straight lines in $L_\diamond$ are {\em streaking artifacts} in the sense of wave front sets. Lemma \ref{lm-fbp} states that such artifacts always exist, namely $\WF(\tilde f_{MA}) \backslash \p D \neq \emptyset$ if the nonlinear function $F$ is quadratic.

 \section{Improved regularity analysis}\label{sec-radon}
 To analyze the singularities produced by  higher order polynomial nonlinearities, we will use a special property of the Radon transform of $\chi_D$ when $D$ is  strictly convex. We start with  Piriou's conormal distributions, see \cite{Pir}, which provides a good motivation. 
 
\begin{definition}  If  $m<-1,$ let $k(m)$ be the non-negative integer such that $-m-2\leq k(m) < -m-1.$ If $\Sigma\subset \Omega$ is a $\CI$ hypersurface, we say that $u\in \ido{}^{m-\frac{n}{4}+\ha}(\Omega,\Sigma)$ if $u\in I^{m-\frac{n}{4}+\ha}(\Omega,\Sigma)$  and vanishes to order $k(m)+1$ at $\Sigma$. 
\end{definition} 

 It is proved in Proposition 2.4 of \cite{SaWa} that 
if  $\Sigma\subset \Omega$ is a $\CI$ hypersurface, $u\in I^{m-\frac{n}{4}+\ha}(\Omega,\Sigma)$ and $m<-1,$ then $u=\mce+ v,$ with  $v \in \ido{}^{m-\frac{n}{4}+\ha}(\Omega,\Sigma)$ and $\mce\in C^\infty.$  If  $v \in \ido{}^{m-\frac{n}{4}+\ha}(\Omega,\Sigma)$ and $\Sigma=\{y_1=0\},$ then $v= y_1^{k(m)} w,$ $w \in I^{m+k(m)-\frac{n}{4}+\ha}(\Omega,\Sigma).$
Now consider $n=2$ and take $v\in \ido{}^{m}(\Omega, \Sigma), m<-1$. Then $v \in y_1^{k(m)}I^{m+k(m)}(\Omega,\Sigma)$.  But since $m+k(m)<-1,$ we get that \begin{gather*}
v^2 \in y_1^{2k(m)} I^{m+k(m)}(\Omega,\Sigma).
\end{gather*}
We can apply Proposition 18.2.3 of \cite{Ho3} 
to conclude that $ 
v^2 \in I^{m-k(m)}(\Omega,\Sigma).
$ 
The argument can be repeated to yield that for $l\in \mbn,$ $v^l \in I^{m - (l-1) k(m)}(\Omega, \Sigma)$. In conclusion, if $k(m)>0$, the conormal distribution becomes more and more regular after self-multiplication. We observe that the vanishing order in Piriou's conormal distribution plays an important role in the argument. 
 
\begin{figure}[t]
\centering
\vspace{-1cm}
\includegraphics[scale = 0.25]{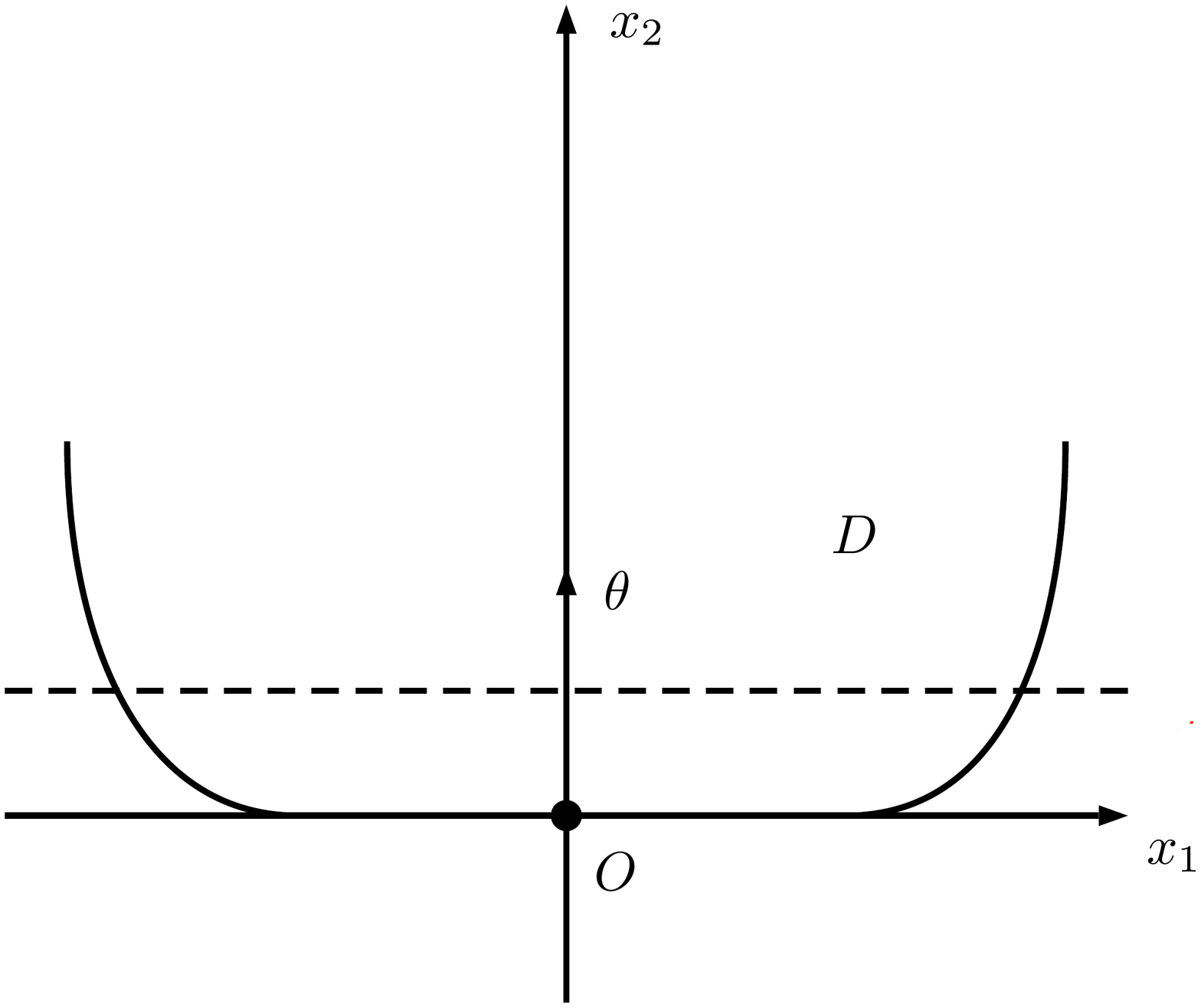}\quad \quad 
\includegraphics[scale = 0.25]{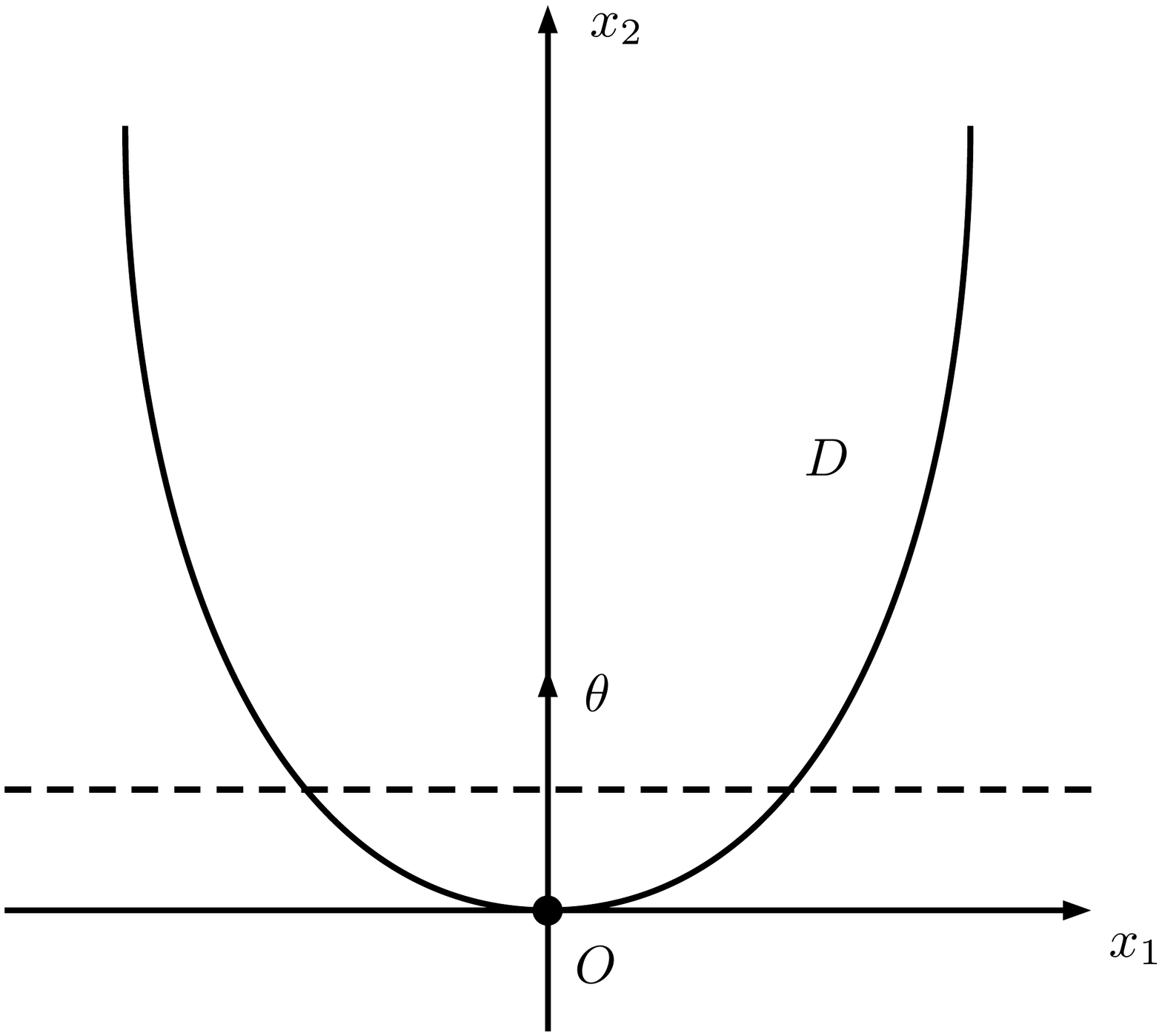}
\vspace{-1cm}
\caption{Regularity of $R\chi_D(s, \theta)$. We consider $\theta$ in the direction of $x_2$-axis. Then $R\chi_D(s, \theta)$ is the integration of $\chi_D(x)$ along the dashed line with distance $s$ to $x_1$-axis. Left figure: $D$ is not strictly convex near $O$.  $R\chi_D(s, \theta)$  has a Heaviside type singularity in $s$. Right figure: $D$ is strictly convex near $O$. Then $R\chi_D(s, \theta)$ behaves like a homogeneous distribution $s_+^{1/2}$ in $s$.}
\label{fig-radon}
\end{figure}

 Now let $U$ be a simply connected bounded domain with smooth boundary $\p U$. As in \eqref{eq-sj}, let $S \defeq \{(s, \phi) \in M : s = x_1\cos\phi + x_2\sin \phi,  x\in \p U \}$ which is a codimension one submanifold of $M.$  We know that $R\chi_U \in I^{-\frac 32}(\mbr^2; S)$ so $m = -3/2$ and $k(m) = 0$. It seems that we do not gain any vanishing  order from the analysis above. However, if $U$ is strictly convex, we show below that it is possible to gain $1/2$ vanishing order.  We remark that for non-strictly convex domain, this is not true. One can construct simple examples to verify it, see Figure \ref{fig-radon}. 

Below, we use $t_+^\alpha, \alpha >-1$ to denote homogeneous distributions so that $t_+^\alpha = t^\alpha$ for $t\geq0$ and $t_+^\alpha = 0$ for $t<0$. See Section 3.2 of \cite{Ho1} for details. The key result of this section is 
\begin{lemma}\label{lm-main}
Suppose $U$ is strictly convex.  For $q_0 \in S$, there exists a neighborhood of $q_0$ and  local coordinates $y = (y_1, y_2)\in \mbr^2$ such that $q_0 = (0, 0)$, $S = \{y_2 = 0\}$ and 
\beq
R\chi_U(y) = 
h(y) y_{2,+}^{1/2} +  y_{2, +} r(y), \quad r\in I^{-3/2}(M;  S)
\eeq 
where $h$ is smooth and positive. 
\end{lemma}

Note that $y_{2, +}\in I^{-3/2}(M; S)$ so the product $y_{2, +}r(y)$ makes sense as product of distributions in $I^{-3/2}(M; S)$ by Lemma \ref{lm-multi}. Also, note that $I^{-3/2}(M; S)\subset L^\infty$, see  \cite{GrU93} so the product also makes sense in $L^\infty.$
\bpf
Recall that $\p U$ is a simple closed strictly convex curve if and only if the curvature $\kappa$ is strictly positive on $\p U$, see Section 2.3 of \cite{Kli}. Here, the curvature is defined in the Frenet frame and is always non-negative. 

For any $q_0 \in S$, we consider a neighborhood and local coordinates $(y_1, y_2)$ such that $q_0 = (0, 0)$ and $S$ is given by $y_2 = 0$. Then we write $s = s(y_1, y_2)$ and $\phi = \phi(y_1, y_2)$ as smooth functions. Consider the Radon transform 
 \beqq\label{eq-calrad}
 \begin{gathered} 
R\chi_D(y_1, y_2)  = \frac{1}{2\pi} \int_{\mbr^2}\int_{\mbr}  e^{i \la (s - x_1 \cos \phi - x_2 \sin\phi)} \chi_D(x) d\la dx_1 dx_2. 
\end{gathered}
 \eeqq 
 For any $\phi$, we let $q$ be a point on $\p U$ such that  $\theta = (\cos \phi, \sin \phi)$  is orthogonal to $T_q(\p U)$. (By the strict convexity of $U$, there are two points with this property.)   We compute the integral in \eqref{eq-calrad} in the Frenet frame at $p$. Thus, we choose local coordinates $(z_1, z_2)$ near $p$ such that $q = (0, 0)$ and $z_1$ axis is tangent to $\p U$ at $q$. Moreover, by selecting the orientation of $\p U$, we can arrange the new coordinate system to  have the same orientation as the original one and $U$ stays in $z_2>0.$  Let $\Phi: \mbr^2\rightarrow \mbr^2$ be the coordinate change so $z = \Phi(x; y)$. Note $\Phi$ depends on $y$ and is smooth in $y.$ Then we write $dx = J(z; y) dz$ where $J$ is the Jacobian factor and $J >0.$ Note that in this new coordinate system, we have $ x_1 \cos \phi + x_2 \sin\phi = s_0  \pm z_2$ 
where the sign is $+$ if $\theta$ is in the direction of positive $z_2$ axis and $-$ otherwise.

Suppose $\p U$  is parameterized by arc-length $\tau$ starting from $q$. Then in the Frenet frame, we have the canonical form of the curve 
 \beq
 z_1(\tau) = \tau - \frac{\kappa^2 \tau^3}{6} + o(\tau^3), \quad z_2(\tau) = \frac{\kappa}{2}\tau^2 + \frac{\kappa' \tau^3}{6} + o(\tau^3)
 \eeq
see Section 1-6 of \cite{Do}. Here, $\kappa, \kappa'$ are the curvature and its $\tau$ derivative at $q.$ As $\kappa>0,$ by using the inverse function theorem, we can take $z_1 \in (-\delta, \delta)$ with $\delta>0$ small as the parameter and express the curve as the graph of a function 
 \beq
 z_2  = \frac{\kappa}{2}z_1^2 + \frac{\kappa' z_1^3}{6}  + o(z_1^3)
 \eeq
 For $z_1 \geq 0$, the function is invertible and we have 
 \beq
 z_1 = (\frac{2}{\kappa})^\ha z_2^\ha + a_3 z_2^{\frac 32} + \cdots, \quad z_2 \geq 0.
 \eeq
 
 Finally, we go back to \eqref{eq-calrad} and get 
  \beq
 \begin{split} 
R\chi_D(y) 
 &= \frac{1}{2\pi} \int_{\mbr}  e^{i \la (s - s_0 \pm z_2)} \int_{\mbr^2} \chi_D(z)  J(z; y) dz d\la \\
 &=   \int_{\mbr}  e^{i \la (s - s_0  \pm z_2)} (\int_0^\delta ( a_1(y, z_2) z_2^\ha + a_3(y, z_2) z_2^{\frac 32} + \cdots)  dz_2) d\la \\
  &=  c_1(y) (s -s_0)^\ha_+ + c_2(y) (s - s_0)^{\frac 32}_+ + \cdots 
\end{split}
 \eeq
 For the second line, we used Taylor expansion of $J(z; y)$ in $z_1.$ Note that $a_1>0$. For the last line, we used the Taylor expansion of $a_n, n = 1, 2, \cdots$ in $z_2$ and integrated in $z_2.$ Note that $c_1(y)>0.$ This completes the proof of the lemma.  
\epf

Next, we consider the situation near $S_\diamond$.
\begin{lemma}\label{lm-main1}
Suppose $D_1, D_2$ are strictly convex.  For $q_0 \in S_1\cap S_2$, there is local coordinates $y = (y_1, y_2)\in \mbr^2$ near $q_0$ such that locally $q_0 = (0, 0),$ $S_1 = \{y_2 = 0\}, S_2 = \{y_1 = 0\}$ and 
\beq
\begin{gathered}
R\chi_{D_1}(y) =  h_1(y) y_{2,+}^{1/2} +  y_{2, +} r_1(y), \quad r_1\in I^{-3/2}(M;  S_1)\\
R\chi_{D_2}(y) =  h_2(y) y_{1,+}^{1/2} +  y_{1, +} r_2(y), \quad r_2 \in I^{-3/2}(M;  S_2)
\end{gathered}
\eeq 
where $h_j, j = 1, 2$ are smooth and positive. \end{lemma}
\bpf
First, we apply Lemma \ref{lm-main1} to find a neighborhood $V_1$ of $q_0$ and coordinates $(z_1, z_2)$ such that $S_1 = \{z_2 = 0\}$. Then we apply Lemma \ref{lm-main1} again to find a neighborhood $V_2$ of $q_0$ and coordinates $(w_1, w_2)$ such that $S_2 = \{w_2 = 0\}$. Because $S_1$ intersects $S_2$ transversally at $q_0$ see \cite{PUW}, we know that the $z_2$ axis is not parallel to $w_2$ axis. Thus, we  can find a new coordinate system $(y_1, y_2)$ with $y_1 = z_2, y_2 = w_2$. Then we write $z_1 = z_1(y), w_1 = w_1(y)$ as smooth functions. Finally, the conclusions follows from Lemma \ref{lm-main1}. 
\epf

The vanishing order is the key for obtaining multiplicative properties similar to Piriou's distributions.  
\begin{lemma}\label{lm-main2}
Let $m, n\in \mbn.$ Under the assumptions of Lemma \ref{lm-main1}, we have 
\begin{enumerate}
\item For $j = 1, 2$, $(R\chi_{D_j})^n(y)  = h_j(y) y_{2,+}^{n/2} +  y_{2, +}^{1+ n/2} r_j(y),$ where $r_j\in I^{-3/2}(M;  S_j)$ and $h_j$ is smooth and positive. 
\item $(R\chi_{D_1})^m(y) (R\chi_{D_2})^n(y) =   h(y) y_{2,+}^{m/2}y_{1,+}^{n/2} + r(y)$,  
where $r$ is a sum of paired Lagrangian distributions such that $r\in I^{-5/2-m/2-n/2}(T_{q_0}^*M\backslash (N^*S_1\cup N^*S_2))$. 
\end{enumerate}
\end{lemma}
\bpf
(1).\ We prove for $n=2$. The general case can be obtained by induction. We can find a proper coordinates as in Lemma \ref{lm-main} and write for $j = 1, 2$
\beq
R\chi_{D_j}(y)  = \tilde h(y) y_{2,+}^{1/2} +  y_{2, +} \tilde r(y), \quad \tilde r\in I^{-3/2}(M;  S_j)
\eeq
Then 
\beq
(R\chi_{D_j})^2 = \tilde h^2(y) y_{2,+} + 2 \tilde h(y) \tilde r(y)  y_{2,+}^{3/2} + y_{2, +}^2 \tilde r^2(y) 
\eeq
Note that $\tilde h \tilde r \in I^{-3/2}(M;  S_j)$.  Also, $\tilde r^2 \in I^{-3/2}(M;  S_j)$ so $y_{2,+}^\ha\tilde r^2 \in I^{-3/2}(M;  S_j)$ as well. 

(2).\ We use the coordinates in Lemma \ref{lm-main1} to get 
\beq
\begin{gathered}
(R\chi_{D_1})^m(y) =  h_1(y) y_{2,+}^{m/2} +  y_{2, +}^{1+ m/2}  r_1(y), \quad r_1\in I^{-3/2}(M;  S_1)\\
(R\chi_{D_2})^n(y) =   h_2(y) y_{1,+}^{n/2} +  y_{1, +}^{1+n/2} r_2(y), \quad r_2 \in I^{-3/2}(M;  S_2)
\end{gathered}
\eeq 
where $h_1, h_2$ are smooth and positive functions. 
Thus
\beq
\begin{gathered}
(R\chi_{D_1})^m(y) (R\chi_{D_2})^n(y)  =   h_1(y) h_2(y) y_{2,+}^{m/2}y_{1,+}^{n/2} +  h_1(y)r_2(y) y_{2,+}^{m/2} y_{1, +}^{1+n/2} \\
+   h_2(y) r_1(y)y_{1,+}^{n/2}y_{2, +}^{1+ m/2}  +   y_{2,+}^{m/2+1}y_{1,+}^{n/2+1} r_1(y)r_2(y).
\end{gathered}
\eeq
Note that $h_1 r_2\in I^{-3/2}(M;  S_2)$. Then by using Proposition 2.4 of \cite{SaWa}, we see that $h_1(y)r_2(y) y_{2,+}^{m/2}\in I^{-1-m/2}(M; S_2)$. By using the proof of Lemma \ref{lm-dij} (here, we need the result for different orders but the proof is the same, see also \cite{GrU93}), we get 
\beqq\label{eq-temp1}
\begin{gathered}
h_1(y)r_2(y) y_{2,+}^{m/2} y_{1, +}^{1+n/2} \in I^{-1-m/2, -3/2-n/2}(T^*_{q_0}M, N^*S_2) \\
+ I^{-2-n/2, -1/2-m/2}(T^*_{q_0}M, N^*S_1).
\end{gathered}
\eeqq
Similarly, we obtain that  
\beqq\label{eq-temp2}
\begin{gathered}
 h_2(y) r_1(y)y_{1,+}^{n/2}y_{2, +}^{1+ m/2}  \in I^{-1-n/2, -3/2-m/2}(T^*_{q_0}M, N^*S_1) \\
 + I^{-2-m/2, -1/2-n/2}(T^*_{q_0}M, N^*S_2),\\
 y_{2,+}^{m/2+1}y_{1,+}^{n/2+1} r_1(y)r_2(y)\in I^{-2-n/2, -3/2-m/2}(T^*_{q_0}M, N^*S_1) \\
 + I^{-2-m/2, -3/2-n/2}(T^*_{q_0}M, N^*S_2).
 \end{gathered}
 \eeqq
 We conclude that away from $N^*S_1\cup N^*S_2$,  terms in \eqref{eq-temp1} and \eqref{eq-temp2} belong  to $I^{-5/2-m/2-n/2}(T_{q_0}^*M)$.  This completes the proof. 
\epf

\section{Determination of the nonlinear term}\label{sec-non} 

\bpf[Proof of Theorem \ref{thm-main}]
Suppose $\tilde F(t) = \sum_{j = 2}^{\tilde J}\tilde a_j t^j$ is another nonlinear polynomial of the form \eqref{eq-F}. Let $\tilde P_{MA}, \tilde f_{CT}$ be the corresponding functions for $\tilde F$. Assume that $f_{CT} = \tilde f_{CT}$. We consider 
\beqq\label{eq-dpma}
P_{MA} - \tilde P_{MA} = \sum_{j = 2}^J (a_j - \tilde a_j)(R\chi_D)^j.
\eeqq
Here, we assumed $J\geq \tilde J$ and we let $\tilde a_j = 0$ for $j>\tilde J.$ 
We claim that for any $q \in S_\diamond$, 
\beqq\label{eq-claim}
\WF(P_{MA} - \tilde P_{MA})\cap (T^*_qM\backslash (N^*S_1\cup N^*S_2)) = \emptyset.
\eeqq 
To see this, we first use Lemma \ref{lm-main2} to conclude that $P_{MA} - \tilde P_{MA}$ is a sum of Lagrangian distributions and paired Lagrangian distributions exactly as those in the proof of Lemma \ref{lm-fbp} but with different orders. Then the symbol calculation in Lemma \ref{lm-fbp}  yields the claim because $f_{CT} - \tilde f_{CT}$ is smooth away from $\p D.$

Next, we  show that $a_j = \tilde a_j$ in \eqref{eq-dpma}. Without loss of generality, we can take $\tilde a_j = 0$ and show $a_j =0.$  We expand and regroup the terms in \eqref{eq-dpma} as following
\beq
\begin{gathered}
\sum_{j = 2}^J  a_j (R\chi_D)^j  = \sum_{j = 2}^J a_j(R\chi_{D_1})^j +  \sum_{j = 2}^J a_j(R\chi_{D_2})^j + \sum_{j = 2}^J a_j A_j\\
\text{ where } A_j = \sum_{m + n  = j, m, n \geq 1} C_{m,n} (R\chi_{D_1})^m (R\chi_{D_2})^n, \quad C_{m, n}>0.
\end{gathered}
\eeq
To determine $a_j,$ we use singularities at $T_q^*M$ for $q\in S_\diamond$ away from $N^*S_1\cup N^*S_2$ so it suffices to look at singularities in $A_j.$ According to Lemma \ref{lm-main2}, we know that  
\beq
A_j = \sum_{m + n  = j, m, n \geq 1} C_{m, n} h_{m, n}(y) y_{2,+}^{m/2}y_{1,+}^{n/2} + r_{m, n}(y),
\eeq
where $C_{m, n}$ and $h_{m,n}$ are both positive. Note that microlocally in $T_q^*M\backslash (N^*S_1\cup N^*S_2),$ we have $C_{m, n} h_{m, n}(y) y_{2,+}^{m/2}y_{1,+}^{n/2} \in I^{-2-m/2-n/2}(T_q^*M)$. We know from Lemma \ref{lm-main2} that $r_{m, n} \in  I^{-5/2-m/2-n/2}(T_q^*M).$ Thus microlocally, $A_j \in I^{-2-j/2}(T_q^*M)$. To see that the principal symbol is non-vanishing, it suffices to find the Fourier transform of $A_{j, 0} = \sum_{m +n =j, m, n \geq 1}C_{m, n} h_{m, n}(0) y_{2,+}^{m/2}y_{1,+}^{n/2}$. Recall from Example 7.1.17 and Section 3.2 of \cite{Ho1} that the Fourier transform of $x_+^a, \re a > -1$ is $\Gamma(a+1) e^{i\pi (a+1)/2}(\xi-i0)^{-a-1}$. 
Thus, the Fourier transform of $A_{j, 0}, j\geq 2$ is given by
\beq
\sum_{m +n =j, m, n \geq 1}C_{m, n} h_{m, n}(0)  \Gamma(\frac m2 +1) \Gamma(\frac n2 + 1)e^{i\pi (j/2+1)/2} (\xi_1-i0)^{-m/2-1} (\xi_2-i0)^{-n/2-1} 
\eeq
Note that $C_{m, n}, h_{m, n}(0)$ and the Gamma function terms are all positive. The exponential factor is the same for all terms in the summation. Thus, the principal symbol of $A_j$ on $T_q^*M\backslash (N^*S_1\cup N^*S_2)$ is non-vanishing. 

 Now we can finish the proof. For $j = 2$, we get that microlocally in $T_q^*M\backslash (N^*S_1\cup N^*S_2)$, $A_2\in I^{-5/2}(T_q^*M)$ and $A_j \in I^{-7/2}(T_q^*M)$ for $j\geq 3$. Because the principal symbol of $A_2$ is non-vanishing, we derive from the claim in the beginning of the proof that $a_2 = 0$. Now we can repeat the argument for $j  = 3, \cdots, J$ to get that all $a_j = 0$. This finishes the proof. 
\epf

 \bpf[Proof of Corollary \ref{cor}]
If $F = 0$, it is easy to see from \eqref{eq-fct} that $f_{CT}\in C^\infty(\mbr^2\backslash \p D).$ If $f_{CT} \in C^\infty$ away from $\p D$, we know from the proof of Theorem \ref{thm-main} that \eqref{eq-claim} holds true for $\tilde P_{MA} = 0$ and $\tilde a_j = 0, j = 1, 2, \cdots, J$. Then the proof of Theorem \ref{thm-main} implies that $F = 0.$ 
\epf

\section{Acknowledgments} 
The author wishes to thank Prof.\ Jin Keun Seo for helpful conversations about the nonlinear nature of the beam-hardening artifacts. This work is  supported by NSF under grant DMS-2205266.
 
\appendix
\section{Composition of FIOs}\label{sec-app}
In this appendix, we verify some technical conditions for the composition of Fourier integral operators in Section \ref{sec-metal}. We recall some definitions from \cite[Section 25.2]{Ho4}. Let $X, Y, Z$ be three manifolds. Let $C_1$ be a homogeneous canonical relation from $T^*Y\backslash 0$ to $T^*X\backslash 0$ and $C_2$ be a homogeneous canonical relation from $T^*Z\backslash 0$ to $T^*Y\backslash 0$, we say that the composition $C_2\circ C_1$ is transversal if $\mcx = C_1\times C_2$ intersects $\mcy = T^*X\times \lap (T^*Y)\times T^*Z$ transversally, that is for any $q$ in the intersection 
\beq
T_q\mcx + T_q \mcy = T_q \mcz \text{ with } \mcz = T^*X\times T^*Y\times T^*Y\times T^*Z
\eeq 
Here,  $\lap(T^*Y)$ denotes the diagonal set of $T^*Y \times T^*Y.$ The composition is proper if the map $\mcx\cap \mcy\rightarrow T^*(X\times Z)\backslash 0$ is proper. If the composition is transversal and proper, then $C =C_2\circ C_1$ is a canonical relation. We have Theorem 25.2.3 of \cite{Ho4} for the composition of FIOs. Actually, we only need the special case of transversal compositions. For the study of composition of Fourier integral operators and conormal distributions, we  take $Z$ to be a point, see the treatment on page 22 of \cite{Ho4}. 

We start with the composition of $C$ and $C_j, j = 1, 2$ in Section \ref{sec-reg}.  
Recall from \eqref{eqcaC} that the canonical relation of $R$ is parametrized as 
\beq 
\begin{gathered}
C  =  \{(x_1 \cos \phi + x_2 \sin \phi, \phi, -\la,  \la(-x_1\sin \phi + x_2 \cos \phi); \\
x_1, x_2,  -\la \cos \phi, - \la \sin \phi):\\
 \la \in \mbr\backslash 0,  \phi \in (-\pi, \pi), x_1, x_2\in \mbr\} \subset T^*M \backslash 0 \times T^*\mbr^2\backslash 0
\end{gathered}
\eeq 
For $C_j = N^*\p D_j\backslash 0, j = 1, 2$, we choose local coordinates $(y_1, y_2)$ near $q\in \p D_j$ such that $q = 0$ and $D_j = \{y_2 = 0\}.$ Then $C_j = \{(y_1, 0, 0, \xi_2): y_1 \in \mbr, \xi_2\in \mbr\backslash 0\}.$ Now we let $\mcz = T^*M\times T^*\mbr^2\times T^*\mbr^2, \mcx = C\times C_j\subset \mcz$ and $\mcy = T^*M\times \lap(T^*\mbr^2) \subset \mcz$.  
Note that $\mcx$ is parametrized by $x_1, x_2, \la, \phi, y_1, \xi_2\in \mbr$ and we write an element of $\mcy$ as $(s, \psi, \alpha, \beta; z_1, z_2, \eta_1, \eta_2,$ $z_1, z_2, \eta_1, \eta_2)$ 
with all variables in $\mbr.$ 
The intersection $\mcx\cap \mcy$ is given by
\beq
\begin{gathered}
s = x_1 \cos \phi + x_2 \sin \phi, \psi = \phi, \alpha = -\la, \beta =  \la(-x_1\sin \phi + x_2 \cos \phi), \\
x_1 = y_1 = z_1, x_2 = y_2 = z_2 = 0, -\la \cos \phi = 0 = \eta_1, -\la \sin \phi = \xi_2 = \eta_2
\end{gathered}
\eeq
which implies that  $\phi = \pm \pi/2$ so 
\beq
\begin{gathered}
\mcx\cap \mcy = \{(0, \pm \pi/2, -\la, \mp\la x_1; x_1, 0, 0, \mp \la, x_1, 0, 0, \mp \la): x_1, \la \in \mbr\}.
\end{gathered}
\eeq
We see that the projection to $T^*M$ is proper. 
Let $q\in \mcx\cap \mcy$. To compute the tangent vector of  $T_q \mcx$, we use the map $\pi: \mbr^6_{(x_1, x_2, \la, \phi, y_1, \xi_2)} \rightarrow \mcx$. So a general tangent vector at $q$ can be obtained by 
\beq
\begin{gathered}
\pi_*(\delta x_1, \delta x_2, \delta \la, \delta \phi, \delta y_1, \delta \xi_2)  \\
= (\delta x_2 - x_1\delta \phi, \delta \phi, -\delta \la, -\la \delta x_1 - x_1 \delta \la, \delta x_1, \delta x_2, \la \delta \phi, -\delta \la, \delta y_1, 0 , 0, \delta \xi_2)
\end{gathered}
\eeq
For the tangent vector of $T_q\mcy$, we use the map $\rho: \mbr^8_{(s, \psi, \alpha, \beta; z_1, z_2, \eta_1, \eta_2)} \rightarrow \mcy$ to get 
\beq
\begin{gathered}
\rho_*(\delta s, \delta \psi, \delta \alpha, \delta \beta; \delta z_1, \delta z_2, \delta \eta_1, \delta \eta_2)\\
 = (\delta s, \delta \psi, \delta \alpha, \delta \beta; \delta z_1, \delta z_2, \delta \eta_1, \delta \eta_2, \delta z_1, \delta z_2, \delta \eta_1, \delta \eta_2)
 \end{gathered}
\eeq 
Now we  conclude that $T_q\mcx + T_q\mcy = T_q \mcz$ by listing $12$ linearly independent tangent vectors which is quite straightforward.

 Next, consider the composition of $C^*$ and $C_0 = T_q^*M\backslash 0$ needed in Lemma \ref{lm-fbp}.  From \eqref{eqcaC}, we get
 \beqq\label{eq-cstar} 
\begin{gathered}
C^*  =  \{(x_1, x_2,  -\la \cos \phi, - \la \sin \phi; x_1 \cos \phi + x_2 \sin \phi, \phi, \\
-\la,  \la(-x_1\sin \phi + x_2 \cos \phi)): 
 \la \in \mbr\backslash 0,  \phi \in (-\pi, \pi), x_1, x_2\in \mbr\} 
\end{gathered}
\eeqq 
We write 
\beq
C_0 = \{(0, 0, \zeta_1, \zeta_2): \zeta_1, \zeta_2 \in \mbr, \zeta_1\zeta_2\neq 0\}
\eeq
Then let $\mcx = C^*\times C_0$, $\mcy = T^*\mbr^2\times \lap(T^*M)$ and $\mcz = T^*\mbr^2\times T^*M\times T^*M.$ The intersection $\mcx\cap \mcy$ is given by
\beq
x_1\cos \phi + x_2\sin \phi = 0, \phi = 0, -\la = \zeta_1, \la(-x_1\sin\phi + x_2\cos \phi) = \zeta_2
\eeq
which implies $x_1 = 0, \phi = 0$ so 
\beq
\mcx\cap \mcy = \{(0, x_2, -\la, 0; 0, 0, -\la, -\la x_2): \la, x_2\in \mbr\}
\eeq
The projection to $T^*\mbr^2$ is proper. Let $q\in \mcx\cap\mcy$. To compute the tangent vector of  $T_q \mcx$, we use the map $\pi: \mbr^6_{(x_1, x_2, \la, \phi, \alpha, \beta)} \rightarrow \mcx$. So 
\beq
\begin{gathered}
\pi_*(\delta x_1, \delta x_2, \delta \la, \delta \phi, \delta \zeta_1, \delta \zeta_2)  \\
= (\delta x_1, \delta x_2, -\delta \la, -\la \delta \phi, \delta x_1+ x_2 \delta \phi, \delta \phi, -\delta \la, x_2\delta \la, \la \delta x_2, 
 0 , 0, \delta \zeta_1, \delta \zeta_2)
\end{gathered}
\eeq
For the tangent vector of $T_q\mcy$, we use $(z_1, z_2, \eta_1, \eta_2; s, \phi, \alpha, \beta; s, \phi, \alpha, \beta)$ for a general element of $\mcy$. Consider  $\rho: \mbr^8_{(z_1, z_2, \eta_1, \eta_2; s, \phi, \alpha, \beta)} \rightarrow \mcy$ and we get 
\beq
\begin{gathered}
\rho_*(\delta z_1, \delta z_2, \delta \eta_1, \delta \eta_2; \delta s, \delta \psi, \delta \alpha, \delta \beta)\\
 = (\delta z_1, \delta z_2, \delta \eta_1, \delta \eta_2, \delta s, \delta \psi, \delta \alpha, \delta \beta, \delta s, \delta \psi, \delta \alpha, \delta \beta)
 \end{gathered}
\eeq 
We can also see that $T_q\mcx + T_q\mcy = T_q \mcz$. 
 
Finally, we consider the composition of $C^*$ and $\tilde C_j = N^*S_j\backslash 0, j = 1, 2$ needed in Lemma \ref{lm-fbp}. In this case, we can find local coordinates $(y_1, y_2)$ near $p\in S_j$ so that $p = 0$ and $S_j = \{y_1 = 0\}$. Thus, 
\beq
\tilde C_j = \{(0, y_2, \eta_1, 0): y_2 \in \mbr, \eta_1\in \mbr\backslash 0 \}
\eeq
Then let $\mcx = C^*\times \tilde C_j$, $\mcy = T^*\mbr^2\times \lap(T^*M)$ and $\mcz = T^*\mbr^2\times T^*M\times T^*M.$ Using \eqref{eq-cstar}, the intersection $\mcx\cap \mcy$ is given by
\beq
x_1\cos \phi + x_2\sin \phi = 0, \phi = y_2, -\la = \eta_1, \la(-x_1\sin\phi + x_2\cos \phi) = 0
\eeq
which implies $x_1 = x_2 = 0$ so 
\beq
\mcx\cap \mcy = \{(0, 0, -\la \cos \phi, - \la \sin\phi, 0, \phi, - \la, 0, 0, \phi, -\la, 0): \la, \phi \in \mbr\}
\eeq
The projection to $T^*\mbr^2$ is proper. To compute the tangent vector of  $T_q \mcx$, we use the map $\pi: \mbr^6_{(x_1, x_2, \la, \phi, y_2, \eta_1)} \rightarrow \mcx$. So 
\beq
\begin{gathered}
\pi_*(\delta x_1, \delta x_2, \delta \la, \delta \phi, \delta \zeta_1, \delta \zeta_2)  \\
= (\delta x_1, \delta x_2, -\delta \la \cos \phi + \la \sin \phi d\phi, -\delta \la \sin \phi - \la \cos \phi \delta \phi; \\
\delta x_1 \cos\phi + \delta x_2 \sin \phi, \delta \phi, -\delta \la, -\la \delta x_1 \sin \phi +  \la \delta x_2 \cos \phi, 
 0 , \delta y_2, -\delta \la, 0)
\end{gathered}
\eeq
For the tangent vector of $T_q\mcy$, we use $(z_1, z_2, \eta_1, \eta_2; s, \phi, \alpha, \beta; s, \phi, \alpha, \beta)$ for a general element of $\mcy$. Let   $\rho: \mbr^8_{(z_1, z_2, \eta_1, \eta_2; s, \phi, \alpha, \beta)} \rightarrow \mcy$ to get
\beq
\begin{gathered}
\rho_*(\delta z_1, \delta z_2, \delta \eta_1, \delta \eta_2; \delta s, \delta \psi, \delta \alpha, \delta \beta)\\
 = (\delta z_1, \delta z_2, \delta \eta_1, \delta \eta_2, \delta s, \delta \psi, \delta \alpha, \delta \beta, \delta s, \delta \psi, \delta \alpha, \delta \beta)
 \end{gathered}
\eeq 
Again, we can find 12 linearly independent vectors to see that $T_q\mcx + T_q\mcy = T_q \mcz$.


\end{document}